\documentclass [11pt]{amsart}
\usepackage {amsmath, amssymb, amscd, graphicx, color, url}
\usepackage[all,knot,arc]{xy}
\newtheorem {theorem}{Theorem}[section]
\newtheorem {definition}[theorem]{Definition}
\newtheorem {lemma}[theorem]{Lemma}
\newtheorem {proposition}[theorem]{Proposition}
\newtheorem {corollary}[theorem]{Corollary}
\newtheorem {conjecture}[theorem]{Conjecture}

\def\zz {{\mathbb{Z}}}

\def\qq {{\mathbb{Q}}}

\def\TT {{\mathcal{T}}}
\def\del {{\partial}}
\def\tt {{\mathfrak{t}}}
\def\ss {{\mathfrak{s}}}
\def\ccc {{\mathcal{C}}}
\def\spc {{\operatorname{Spin^c}}}
\def\sgn {{\operatorname{sgn}}}
\def\rk {{\operatorname{rank}}}
\def\fin {{\hfill \square}}

\newenvironment{narrow}[2]{%
 \begin{list}{}{%
  \setlength{\topsep}{0pt}%
  \setlength{\leftmargin}{#1}%
  \setlength{\rightmargin}{#2}%
  \setlength{\listparindent}{\parindent}%
  \setlength{\itemindent}{\parindent}%
  \setlength{\parsep}{\parskip}%
 }%
\item[]}{\end{list}}

\newif\ifpic
\picfalse    % no pictures
\pictrue   % with pictures

\begin{document}

\title
[A concordance invariant from Floer homology]{A concordance invariant from the Floer homology of
double branched covers}

\author [Ciprian Manolescu]{Ciprian Manolescu}
\thanks {CM was supported by a Clay Research Fellowship.}
\address {Department of Mathematics, Columbia University\\ New York, NY 10027}
\email {cm@math.columbia.edu}

\author [Brendan Owens]{Brendan Owens}
\address {Department of Mathematics, Cornell University\\ Ithaca, NY 14853}
\email {owensb@math.cornell.edu}

\begin {abstract} Ozsv\'ath and Szab\'o defined an analog of the Fr\o yshov invariant in the form of a 
correction term for the grading in Heegaard Floer homology. Applying this to the double cover of the 
3-sphere branched over a knot $K$, we obtain an invariant $\delta$ of knot concordance. We show that 
$\delta$ is determined by the signature for alternating knots and knots with up to nine crossings, and 
conjecture a similar relation for all H-thin knots. We also use $\delta$ to prove that for all knots $K$ 
with $\tau(K) > 0,$ the positive untwisted double of $K$ is not smoothly slice. \end {abstract}

\maketitle

\section {Introduction}

In \cite{OS2}, Ozsv\'ath and Szab\'o associated an invariant  $d(Y, \tt) \in \qq$ to every rational
homology three-sphere $Y$ endowed with a $\spc$ structure. In this paper we study the knot invariant
$$\delta(K) = 2d(\Sigma(K), \tt_0),$$ where $\Sigma(K)$ is the double cover of $S^3$ branched over the
knot $K,$ and $\tt_0$ is the $\spc$ structure induced by the unique Spin structure on $\Sigma(K).$
In \cite{OS2} it is shown that $d$ induces a group homomorphism from the three-dimensional $\spc$
homology bordism group to $\qq.$ An immediate consequence of this fact and the basic properties of
$d$ is the following:

\begin {theorem}
\label {theo}
The invariant $\delta(K)$ descends to give a surjective group homomorphism $\delta: \ccc \to \zz,$
where $\ccc$ is the smooth concordance group of knots in $S^3.$
\end {theorem}

It is interesting to compare $\delta$ to three other homomorphisms from $\ccc$ to the integers. The
first is the classical knot signature $\sigma,$ which we normalize to $\sigma' = - \sigma/2.$ The
second is the invariant $\tau$ defined using the knot Floer homology of Ozsv\'ath-Szab\'o and
Rasmussen \cite{OS6}, \cite{R2}. The third is Rasmussen's invariant $s$ coming from Khovanov homology
\cite{R3}, which we normalize to $s' = -s/2.$ For alternating knots it is known that $\tau = s' =
\sigma'.$ We show that a similar result holds for $\delta:$

\begin {theorem}
\label {alt}
If the knot $K$ is alternating, then $\delta(K) = \sigma'(K).$
\end {theorem}

The similarities between the four invariants hold for a much larger class of knots. Indeed, in \cite{R3} 
Rasmussen conjectured that $s' = \tau$ for all knots. On the other hand, there are several known examples 
where $s' = \tau \neq \sigma'.$ Following \cite{Kh}, we call a knot H-thin (homologically thin) if its 
Khovanov homology is supported on two adjacent diagonals, and H-thick otherwise. Alternating knots are 
H-thin by the work of Lee \cite{Le}, as are most knots up to ten crossings \cite{BN}, \cite{Kh}. For all 
H-thin knots for which $s'$ was computed, it turned out to be equal to $\sigma'.$

The invariant $\delta$ can be computed algorithmically for Montesinos and torus knots, using the fact that 
their double branched covers are Seifert fibrations. We performed this computation for all Montesinos knots 
in a certain range, and found that $\delta = \sigma'$ for most of them. Interestingly, all the exceptions 
were H-thick. We also got $\delta= \sigma'$ for many H-thick examples, such as $8_{19}$ and $9_{42}.$

In another direction, we used the method of Ozsv\'ath and Szab\'o from \cite{OS8} to compute $\delta$ for  
several H-thin knots of a special form. Together with the calculations for alternating, 
Montesinos and torus knots, this is enough to cover all the knots with up to nine crossings. We find:

\begin {theorem}
\label {nine}
If the knot $K$ admits a diagram with nine or fewer crossings, then $\delta(K) = \sigma'(K).$
\end {theorem}

We also managed to calculate $\delta$ for all but eight of the ten-crossing knots. Among these we only found 
two examples with $\delta$ different from $\sigma',$ the H-thick knots $10_{139}$ and $10_{145}.$
Based on Theorem~\ref{alt} and all our computations, we make the following:

\begin {conjecture}
\label {conj}
For any H-thin knot $K, \ \delta(K) = \sigma'(K).$
\end {conjecture}

On the other hand, an important difference between $\delta$ and the other three invariants $\sigma',
\tau,$ and $s'$ is that the absolute value of $\delta$ does not provide a lower bound for the slice
genus of a knot. For example, the slice genus of the knot $10_{145}$ is two, but $\delta(10_{145}) =
-3.$ However, since $\delta$ is a concordance invariant, it can still be used as an obstruction to
sliceness.

We can also say something about $\delta$ for Whitehead doubles. We denote by $Wh(K)$ the untwisted
double of $K$ with a positive clasp. We show:

\begin {theorem}
\label {double}
For any knot $K$ we have $\delta(Wh(K)) \leq 0,$ and the inequality is strict if $\tau(K) >
0.$ If $K$ is alternating, then $\delta(Wh(K)) = -4\max \{\tau(K), 0\}.$
\end {theorem}

It is well-known that the Alexander polynomial of $Wh(K)$ is
always $1,$ and consequently $Wh(K)$ is topologically slice \cite{FK}. Our interest lies in
the following conjecture, which appears as Problem 1.38 in Kirby's list \cite{Ki}.

\begin {conjecture}
$Wh(K)$ is (smoothly) slice if and only if $K$ is slice.
\end {conjecture}
A quick corollary of Theorem~\ref{double} is a result in the direction of this conjecture:
\begin {corollary}
\label {coro}
If the knot $K$ has $\tau(K) > 0,$ then $Wh(K)$ is not slice.
\end {corollary}

Previously, Rudolph proved in \cite{Ru} that $Wh(K)$ is not slice whenever the Thurston-Bennequin 
invariant $TB(K)$ of the knot $K$ is nonnegative. (See also \cite{AM}, \cite{L} for different 
proofs.) Since $TB(K) \leq 2\tau(K) - 1$ by the work of Plamenevskaya \cite{P}, Rudolph's result is a 
consequence of ours. Furthermore, Corollary~\ref{coro} can also be applied to knots such as $6_2, 
7_6, 8_{11}$ and the mirrors of $8_4, 8_{10}, 8_{16},$ which according to the knot table in \cite{Lt} 
have $\tau > 0$ but $TB < 0.$

As suggested to us by Jacob Rasmussen, Theorem~\ref{double} can be used to produce examples of knots for 
which $\delta$ is a nontrivial obstruction to sliceness, while the other three invariants are not. 
Indeed, we have:

\begin {corollary}
\label {coronew}
Let $K_1=T(2,2m+1)$ and $K_2=T(2,2n+1)$ be two positive torus knots with $m, n \geq 1, \ m\neq n.$ 
Then the connected sum $Wh(K_1) \# (-Wh(K_2))$ has $\sigma' = \tau = s' =0$ but $\delta \neq 0.$
\end {corollary}

Let us denote by $\ccc_{ts}$ the smooth concordance group of topologically slice knots. It was shown by 
Livingston \cite{L} that $\tau: \ccc_{ts} \to \zz$ is a surjective homomorphism, and therefore $\ccc_{ts}$ 
has a $\zz$ summand. Using the pair $(\delta, \tau),$ we show:

\begin {corollary}
\label {coro2}
The group $\ccc_{ts}$ has a $\zz \oplus \zz$ summand.
\end {corollary}

\medskip
\noindent \textbf{Acknowledgements.} We wish to thank Matt Hedden and Jacob Rasmussen for some very helpful 
discussions, and Alexander Shumakovitch and Sa\v{s}o Strle for computer help. We are also grateful to Chuck 
Livingston for suggesting to us Corollary~\ref{coro2} above, and for correcting a misstatement in a previous 
version of the paper.

\section {General facts}

In this section we prove Theorem~\ref{theo}. First, we claim that
$\delta(K)$ is an integer for any knot $K.$ Note that the double
branched cover $Y= \Sigma(K)$ is a rational homology three-sphere
with $|H_1(Y)| = \det(K)$ an odd integer. Any $Y$ with $|H_1(Y)|$
odd admits a unique Spin structure $\tt_0.$ This can be
distinguished in the set of all $\spc$ structures by the
requirement that $c_1(\tt_0) = 0 \in H^2(Y; \zz).$ Let $X$ be a
four-manifold equipped with a $\spc$ structure $\ss$ such that
$\del X = Y$ and $\ss|_{Y} = \tt_0.$ Then $k = c_1(\ss)$ is an
element in the kernel of the map $H^2(X; \zz) \to H^2(Y; \zz).$
According to \cite{OS2}, we have
$$  \delta(K) = 2d(Y, \tt_0) \equiv \frac{k^2 - \sgn(X)}{2} \,\pmod4.$$
Here $\sgn(X)$ denotes the signature of $X,$ and $k$ is a characteristic element for the intersection
form on $X.$

We choose $X$ to be the double cover of $B^4$ branched along the pushoff of a Seifert surface for $K$. It is
shown in \cite{Kauf} that $X$ is a spin manifold (so that $k=0$ is characteristic) with $\sgn(X)=\sigma(K)$.
It follows that $\delta(K)$ is an integer and, furthermore,
\begin {equation}
\label {congruence}
\delta(K)\equiv \sigma'(K)\,\pmod 4.
\end {equation}

Next, assume that $K$ is a slice knot, i.e. it sits on the boundary $S^3$ of $B^4$ and bounds a smooth disk $D \subset
B^4.$ Let $X$ be the double cover of $B^4$ branched over $D.$ Then $X$ is a rational homology four-ball with boundary $Y =
\Sigma(K).$ Since $d$ is an invariant of $\spc$ rational homology bordism \cite{OS2}, it follows that $d(Y, \tt_0) = 0$
and therefore $\delta(K) = 0.$

The additivity property for connected sums $\delta(K_1 \# K_2) = \delta(K_1) + \delta(K_2)$ is a
consequence of the additivity of $d$ (Theorem~4.3 in \cite{OS2}). Also, Proposition 4.2 in \cite{OS2}
implies that $\delta(-K) = - \delta(K),$ where by $-K$ we denote the mirror of $K.$ If $K_1$ and $K_2$
are two cobordant knots, then $K_1 \# (-K_2)$ is slice, hence $\delta(K_1) = \delta(K_2).$ This shows
that $\delta$ is a well-defined group homomorphism from the concordance group $\ccc$ to $\zz.$
Surjectivity is an easy consequence of Theorem~\ref{alt} below. For example, $\delta(T(3,2)) = 1$ for
the right-handed trefoil.

\section {Alternating knots}
\label{sec:alt}

This section contains the proof of Theorem~\ref{alt}. The main input comes from the paper \cite{OS7}, where
Ozsv\'ath and Szab\'o calculated the Heegaard Floer homology for the double branched covers of alternating
knots.

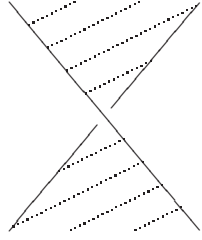
\begin{figure}[htbp] %% colouring
\begin{center}
\ifpic
\leavevmode
\begin{xy}
0;/r6pc/:
% Name some points
(0,0)*{}="1";
(1,0)*{}="2";
(1,1.2)*{}="3";
(0,1.2)*{}="4";
"2";"4" **\crv{}; \POS?(.5)*{\hole}="x"; % Draw overcrossing, mark a point halfway
"1";"x" **\crv{}; % draw lower piece of undercrossing
"x";"3" **\crv{}; % draw upper piece of undercrossing
% Now add shading (dotted lines)
(0.1,1.08);(0.34,1.2) **\dir{.};
(0.2,0.96);(0.68,1.2) **\dir{.};
(0.3,0.84);(0.98,1.18) **\dir{.};
(0.4,0.72);(0.74,0.89) **\dir{.};
(0.6,0.48);(0.26,0.31) **\dir{.};
(0.7,0.36);(0.02,0.02) **\dir{.};
(0.8,0.24);(0.32,0) **\dir{.};
(0.9,0.12);(0.66,0) **\dir{.};
\end{xy}
\else \vskip 5cm \fi
\begin{narrow}{0.3in}{0.3in}
\caption{
\bf{Coloring convention for alternating knot diagrams.}}
\label{coloring}
\end{narrow}
\end{center}
\end{figure}

Let $K$ be a knot with a regular, alternating projection. The projection splits the two-sphere into several regions,
which we color black and white in chessboard fashion. Our coloring convention is that at each crossing, the white
regions should be to the left of the overpass (see Figure~\ref{coloring}). We form a graph $\Gamma$ as follows.  We
form the set of vertices $V(\Gamma)$ by assigning a vertex to each white region. We index the vertices by setting
$V(\Gamma) = \{v_0,v_1, \dots, v_m\},$ so that the number of white regions is $m+1.$ Let us denote by $R_i$ the
white
region corresponding to $v_i.$ We draw an edge between the vertices $i, j \in V(\Gamma)$ for every crossing at which
the regions $R_i$ and $R_j$ come together.  We denote by $e_{ij}$ the number of edges joining the vertices $v_i$ and
$v_j.$ We can assume that the diagram contains no reducible crossings, and therefore $e_{ii} = 0$ for all $i.$

\begin{figure}[htbp]
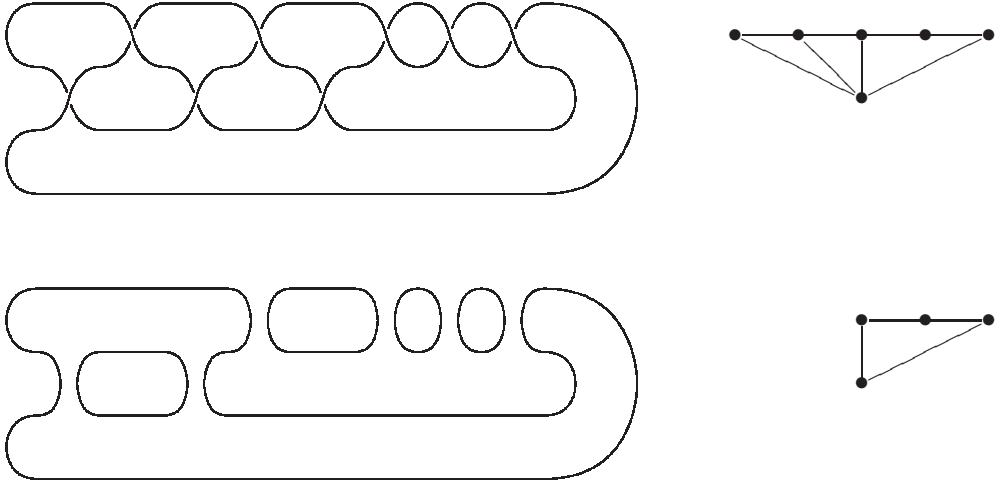
  %%% knot 8_13, oriented resolution, graphs
\begin{center}
\ifpic
\leavevmode
\xygraph{
!{0;/r2.0pc/:}
%%% first the knot
%open
!{\hcap[-1]}[dd]
!{\hcap[-1]}[uu]
%negativemiddle
!{\xcaph[1]@(0)}[ld]
!{\htwistneg}[ldd]
!{\xcaph[1]@(0)}[uuu]
%positivetop
!{\htwist}[ldd]
!{\xcaph[1]@(0)}[ld]
!{\xcaph[1]@(0)}[uuu]%}
%negativemiddle
!{\xcaph[1]@(0)}[ld]
!{\htwistneg}[ldd]
!{\xcaph[1]@(0)}[uuu]
%positivetop
!{\htwist}[ldd]
!{\xcaph[1]@(0)}[ld]
!{\xcaph[1]@(0)}[uuu]%}
%negativemiddle
!{\xcaph[1]@(0)}[ld]
!{\htwistneg}[ldd]
!{\xcaph[1]@(0)}[uuu]
%positivetop
!{\htwist}[ldd]
!{\xcaph[1]@(0)}[ld]
!{\xcaph[1]@(0)}[uuu]%}
%positivetop
!{\htwist}[ldd]
!{\xcaph[1]@(0)}[ld]
!{\xcaph[1]@(0)}[uuu]%}
%positivetop
!{\htwist}[ldd]
!{\xcaph[1]@(0)}[ld]
!{\xcaph[1]@(0)}[uuu]%}
%close(even length contfrac)
!{\hcap[3]}[d]
!{\hcap}
%%% now draw the white graph
%put some points and name them
[r(3)u(0.5)]!{*{\bullet}}="v1"
[r]!{*{\bullet}}="v2"
[r]!{*{\bullet}}="v3"
[r]!{*{\bullet}}="v4"
[r]!{*{\bullet}}="v5"
[l(2)d]!{*{\bullet}}="v6"
%put edges
"v1"-"v5"
"v1"-"v6"
"v2"-"v6"
"v3"-"v6"
"v5"-"v6"[l(13)d(3)]
%%% now the oriented resoution
%open
!{\hcap[-1]}[dd]
!{\hcap[-1]}[uu]
%uncrossmiddle
!{\xcaph[1]@(0)}[ld]
!{\huncross}[ldd]
!{\xcaph[1]@(0)}[uuu]
%fourlines
!{\xcaph[1]@(0)}[ld]
!{\xcaph[1]@(0)}[ld]
!{\xcaph[1]@(0)}[ld]
!{\xcaph[1]@(0)}[uuu]%}
%uncrossmiddle
!{\xcaph[1]@(0)}[ld]
!{\huncross}[ldd]
!{\xcaph[1]@(0)}[uuu]
%uncrosstop
!{\huncross}[ldd]
!{\xcaph[1]@(0)}[ld]
!{\xcaph[1]@(0)}[uuu]
%fourlines
!{\xcaph[1]@(0)}[ld]
!{\xcaph[1]@(0)}[ld]
!{\xcaph[1]@(0)}[ld]
!{\xcaph[1]@(0)}[uuu]%}
%uncrosstop
!{\huncross}[ldd]
!{\xcaph[1]@(0)}[ld]
!{\xcaph[1]@(0)}[uuu]
%uncrosstop
!{\huncross}[ldd]
!{\xcaph[1]@(0)}[ld]
!{\xcaph[1]@(0)}[uuu]
%uncrosstop
!{\huncross}[ldd]
!{\xcaph[1]@(0)}[ld]
!{\xcaph[1]@(0)}[uuu]
%close(even length contfrac)
!{\hcap[3]}[d]
!{\hcap}
%%% draw the resolved graph
%put some points and name them
[r(5)u(0.5)]!{*{\bullet}}="v3"
[r]!{*{\bullet}}="v4"
[r]!{*{\bullet}}="v5"
[l(2)d]!{*{\bullet}}="v6"
%put edges
"v3"-"v5"
"v3"-"v6"
"v5"-"v6"
}
\else \vskip 5cm \fi
\begin{narrow}{0.3in}{0.3in}
\caption{
{\bf{The knot $8_{13}$.}}
The top row contains the knot diagram and the white graph $\Gamma$, and the second row
contains the oriented resolution and the graph $\tilde \Gamma$.}
\label{whitegraph}
\end{narrow}
\end{center}
\end{figure}

Choose an orientation on $K,$ and consider the oriented resolution of the knot projection. Note that at every positive
crossing the two white regions get joined, while at every negative crossing they are separated. Now form a new graph
$\tilde \Gamma$ by taking its set of vertices $V(\tilde \Gamma)$ to be the set of white regions in the resolved
diagram, and by assigning exactly one edge whenever two white domains are adjacent at some negative crossing in the
original knot projection. In other words, we look at the natural projection map $\pi : V(\Gamma) \to V(\tilde \Gamma).$
Given $x\ne y \in V(\tilde \Gamma),$ there is at most one edge between $x$ and $y.$ Such an edge exists if and only if
there are $v_i \in \pi^{-1}(x), v_j\in \pi^{-1}(y)$ such that $v_i$ and $v_j$ are joined by
at least one edge in $G.$ This is shown in Figure~\ref{whitegraph} for the two-bridge knot $8_{13}$.

A circuit in $\tilde \Gamma$ corresponds to a (black) region in the original knot diagram
with only negative crossings on its boundary.  Since the knot is alternating, it follows that
all circuits contain an even number of edges.
This means that $\tilde \Gamma$ can be made into a bipartite graph as follows. The distance
between two vertices on a graph is defined to be the minimal number of edges in a path from one vertex to the other.
We define a map $\tilde \eta: V(\tilde \Gamma) \to \{0,1\}$ by setting $\tilde \eta (x) = $ the parity of the
distance between $\pi(0)$ and $x.$ This produces a partition of $V(\tilde \Gamma)$ into two classes such that
$\tilde \eta(x) \neq \tilde \eta(y)$ whenever $x$ and $y$ are joined by an edge. We can pull back $\tilde \eta$ to a
partition of $V(\Gamma)$ described by the map $\eta = \tilde \eta \circ \pi: V(\Gamma) \to \{0, 1\}.$

\begin {lemma}
\label {cross}
An edge in $\Gamma$ between $v_i$ and $v_j$ corresponds to a positive crossing if and only if $\eta(v_i)
= \eta(v_j).$
\end {lemma}

\noindent {\bf Proof.} If two domains $R_i$ and $R_j$ meet at a positive crossing, then $\pi(v_i) = \pi(v_j),$
which implies $\eta(v_i)= \eta(v_j).$

If $R_i$ and $R_j$ meet at a negative crossing $c$, we claim that any other crossing $c'$ where $R_i$ and $R_j$ meet
is also negative. Indeed, let us draw a loop $\gamma$ on the two-sphere by going from $c$ to $c'$ inside $R_i$ and
then going from $c'$ to $c$ in $R_j.$ Let $D$ one of the two disks bounded by $\gamma.$ Since $c$ is negative, the
underpass and the overpass at $c$ must be going either both into $D$ or both out of $D.$ It follows that the same is
true for the underpass and the overpass at $c',$ which means that $c'$ is negative. This proves our claim.

Therefore, if $R_i$ and $R_j$ meet at a negative crossing, the corresponding $v_i$ and $v_j$ must project to
different vertices $x = \pi(v_i), y =\pi(v_j) \in V(\tilde \Gamma),$ where $x$ and $y$ are joined by an edge in
$\tilde \Gamma.$ We observed above that in this situation $\tilde \eta(x) \neq \tilde \eta(y). \hfill \fin$
\medskip

The Goeritz matrix $G = (g_{ij})$ of the alternating projection is the $m \times m$ symmetric
matrix with entries:
$$g_{ij} = \begin {cases} -\sum_{k=0}^m e_{ik} & \text{if }i=j; \\
  e_{ij} & \text{if } i \neq j.
\end {cases}$$

The quadratic form associated to $G$ is negative-definite. Indeed, if we have a vector $w = (w_i), i=1 \dots
m,$ then

\begin {equation}
\label {goeritz}
 w^{T}Gw = \sum_{i=1}^m \sum_{j=1}^m g_{ij} w_i w_j = -\sum_{i=1}^m e_{0i} w_i^2 - \sum_{1\leq i < j \leq m}
e_{ij} (w_i - w_j)^2.
\end {equation}

\begin {definition}
\label{char}
Given an $m \times m$ symmetric matrix $M = (m_{ij})$
with integer entries, a vector $w = (w_i), i=1\dots m,$ with $w_i
\in \zz$ is called {\bf characteristic} for $M$ if $(Mw)_i \equiv
m_{ii} \pmod 2$ for all $i.$
\end {definition}

\begin {lemma}
\label {charmod}
If $w = (w_i)$ is a characteristic vector for the Goeritz matrix $G$ of an alternating projection, then $w_i \equiv
\eta(v_i) \pmod 2.$
\end {lemma}

\noindent {\bf Proof.} The determinant of $G$ is the same as the determinant of $K$ (c.f. \cite[Corollary 13.29]{BZ}),
which is an odd number. This means that the system of equations
\begin {equation}
\label {congr}
 \sum_{j=1}^m g_{ij}w_j = g_{ii} \ \ \ (i=1, \dots, m)
\end {equation}
has a unique solution over the field $\zz/2.$ Therefore, it suffices to show that $w_i = \eta(v_i)$ satisfy
(\ref{congr}). Using the fact that $\eta(v_0) =0,$ we can rewrite this condition as:
$$ \sum_{j=0}^m e_{ij} (\eta(v_j) - \eta(v_i) +1) \equiv 0 \pmod 2\  \text{ for } i=1, \dots, m.$$

The left hand side of this congruence is the number of positive crossings around the white region $R_i.$ The
boundary of $R_i$ is a polygon with an orientation on each edge. A corner of $R_i$ corresponds to a positive
crossing if and only if the orientations of the two incident edges do not match up at that corner. The number of
those corners has to be even. $\hfill \fin$

\medskip

Gordon and Litherland \cite[Theorem 6]{GL} proved a formula for the signature of a knot in terms of the Goeritz
matrix. In the case of an alternating projection (with our coloring convention), their formula reads:
\begin {equation}
\label {signature}
\sigma(K) = - m + \mu,
\end {equation}
where $\mu$ is the number of negative crossings.

On the other hand, the formula for the correction term $d(\Sigma(K), \tt_0)$ in \cite[Theorem 3.4]{OS7} implies:
\begin {equation}
\label {minm}
2\delta(K) = m+ \max_{w} \{w^T G w\},
\end {equation}
 where the
maximum is taken over characteristic elements $w.$ From equation (\ref{goeritz}) we see that the maximum can be
attained by choosing all the $w_i$ to be either $0$ or $1.$ Using the description of the characteristic vectors in
Lemma~\ref{charmod}, it follows that we should take $w_i = \eta(v_i)$ for all $i.$ Equation (\ref{goeritz}) shows
that the maximum appearing in (\ref{minm}) is minus the total number of edges between vertices $v_i$ and $v_j$ with
$\eta(v_i) \neq \eta(v_j).$ By Lemma~\ref{cross}, this is exactly $-\mu.$ Comparing (\ref{signature}) with
(\ref{minm}), we get that $\delta(K) = -\sigma(K)/2.$ This completes the proof of Theorem~\ref{alt}.

\section {Montesinos and torus knots}

\subsection {Montesinos knots}

For a detailed exposition of the properties of Montesinos knots
and links we refer to \cite{BZ}. In the following definition,
assume that $e$ is any integer and $(\alpha_i, \beta_i)$ are
coprime pairs of integers with $\alpha_i > \beta_i \geq 1$.

\begin{definition}
\label{def:mont} A Montesinos link
$M(e;(\alpha_1,\beta_1),(\alpha_2,\beta_2),\ldots,(\alpha_r,\beta_r))$
is a link which has a projection as shown in Figure
\ref{fig:mont}(a).  There are $e$ half-twists on the left side.  A
box \framebox{$\alpha,\beta$} represents a {\em rational tangle of
slope} $\alpha/\beta$:  given a continued fraction expansion
$$\frac{\alpha}{\beta}=[a_1,a_2,\ldots,a_m]:=a_1-\frac{1}{a_2-\raisebox{-3mm}{$\ddots$
\raisebox{-2mm}{${-\frac{1}{\displaystyle{a_m}}}$}}}\, ,$$ the
rational tangle of slope $\alpha/\beta$ consists of the four
string braid
$\sigma_2^{a_1}\sigma_1^{a_2}\sigma_2^{a_3}\sigma_1^{a_4}\ldots\sigma_i^{a_m}$,
which is then closed on the right as in Figure \ref{fig:mont}(b)
if $m$ is odd or (c) if $m$ is even.
\end{definition}

\begin{figure}[htbp]
\begin{center}
\ifpic \leavevmode \xygraph{ !{0;/r1.0pc/:}
%first label
!{*{(a)}}[rd]
%top of montesinos
!{\vcap[4]}[r] !{\vcap[2]}[l]
%left half of montesinos
!{\xcapv@(0)}[ur] !{\xcapv@(0)}[l] !{\xcapv@(0)}[ur]
!{\xcapv@(0)}[l] !{\vtwist} !{\vtwist} !{\vtwist}
!{\xcapv@(0)}[ur] !{\xcapv@(0)}[l] !{\xcapv@(0)}[ur]
!{\xcapv@(0)}[l]
%bottom of montesinos
!{\vcap[-4]}[r] !{\vcap[-2]}[r(2.5)u(6.3)]
%right half of montesinos
!{*+[F]{\alpha_1,\beta_1}}[l(0.5)d(0.75)] !{\xcapv[0.5]@(0)}[ur]
!{\xcapv[0.5]@(0)}[l(0.5)d(0.2)]
!{*+[F]{\alpha_2,\beta_2}}[l(0.5)d(0.75)] !{\xcapv[0.5]@(0)}[ur]
!{\xcapv[0.5]@(0)}[l(0.5)u(0.1)] !{*{.}}[d(0.2)] !{*{.}}[d(0.2)]
!{*{.}}[d(0.35)l(0.5)] !{\xcapv[0.5]@(0)}[ur]
!{\xcapv[0.5]@(0)}[l(0.5)d(0.2)]
!{*+[F]{\alpha_r,\beta_r}}[r(5)u(7.25)]
%now do first version of tangle
%second label
!{*{(b)}}[r]
%curving up
!{\xbendd[-2]@(0)}[u] !{\xbendd-@(0)}[dl] !{\xcaph[1]@(0)}[l(2)d]
%curving down
!{\xbendu[-2]@(0)}[d] !{\xbendu-@(0)}[ul]
!{\xcaph[1]@(0)}[l(0.5)u(2)]
%notwists
!{\xcaph[1]@(0)}[ld] !{\xcaph[1]@(0)}[ld] !{\xcaph[1]@(0)}[ld]
!{\xcaph[1]@(0)}[uuu]
%positivemiddle
!{\xcaph[1]@(0)}[ld] !{\htwist}[ldd] !{\xcaph[1]@(0)}[uuu]
%positivemiddle (with label)
!{\xcaph[1]@(0)}[ld] !{\htwist|{a_1}}[ldd] !{\xcaph[1]@(0)}[uuu]
%positivemiddle
!{\xcaph[1]@(0)}[ld] !{\htwist}[ldd] !{\xcaph[1]@(0)}[uuu]
%negativetop (with label clumsily shifted to right)
!{\htwistneg>{\,\,\,\,\,\,\,\,\,\,\,\,\,\,\,\,\,\,\,\,a_2}}[ldd]
!{\xcaph[1]@(0)}[ld]
!{\xcaph[1]@(0)}[uuu]%}
%negativetop
!{\htwistneg}[ldd] !{\xcaph[1]@(0)}[ld]
!{\xcaph[1]@(0)}[uuu]%}
%positivemiddle (with label)
!{\xcaph[1]@(0)}[ld] !{\htwist|{a_3}}[ldd] !{\xcaph[1]@(0)}[uuu]
%close(odd length contfrac)
!{\hcap}[dd] !{\hcap}[uu] [r(3)u]
%then second version of tangle
%third label
!{*{(c)}}[r]
%curving up
!{\xbendd[-2]@(0)}[u] !{\xbendd-@(0)}[dl] !{\xcaph[1]@(0)}[l(2)d]
%curving down
!{\xbendu[-2]@(0)}[d] !{\xbendu-@(0)}[ul]
!{\xcaph[1]@(0)}[l(0.5)u(2)]
%notwists
!{\xcaph[1]@(0)}[ld] !{\xcaph[1]@(0)}[ld] !{\xcaph[1]@(0)}[ld]
!{\xcaph[1]@(0)}[uuu]
%positivemiddle
!{\xcaph[1]@(0)}[ld] !{\htwist}[ldd] !{\xcaph[1]@(0)}[uuu]
%positivemiddle (with label)
!{\xcaph[1]@(0)}[ld] !{\htwist|{a_1}}[ldd] !{\xcaph[1]@(0)}[uuu]
%positivemiddle
!{\xcaph[1]@(0)}[ld] !{\htwist}[ldd] !{\xcaph[1]@(0)}[uuu]
%negativetop
!{\htwistneg}[ldd] !{\xcaph[1]@(0)}[ld] !{\xcaph[1]@(0)}[uuu]
%negativetop (with label)
!{\htwistneg>{a_2}}[ldd] !{\xcaph[1]@(0)}[ld]
!{\xcaph[1]@(0)}[uuu]
%negativetop
!{\htwistneg}[ldd] !{\xcaph[1]@(0)}[ld] !{\xcaph[1]@(0)}[uuu]
%close(even length contfrac)
!{\hcap[3]}[d] !{\hcap} } \else \vskip 5cm \fi
\begin{narrow}{0.4in}{0.2in}
\caption{ {\bf{Montesinos links and rational tangles.}} Note that
$e=3$ in (a).  Also (b) and (c) are both representations of the
rational tangle of slope 10/3:
$$10/3=[3,-2,1]=[3,-3]$$ (and one can switch between (b) and (c) by simply
moving the last crossing).} \label{fig:mont}
\end{narrow}
\end{center}
\end{figure}

Let $K$ be the link
$M(e;(\alpha_1,\beta_1),(\alpha_2,\beta_2),\ldots,(\alpha_r,\beta_r))$.
Note that $K$ is alternating if $e\notin\{1,2,\dots,r-1\}$: for an
alternating projection in the case $e\le0$, take a continued
fraction expansion
$$\frac{\alpha_i}{\beta_i}=[a^i_1,a^i_2,\ldots,a^i_{m_i}]$$
with $a^i_1,a^i_3,\dots$ positive and $a^i_2,a^i_4,\dots$
negative, for each $i$ (as in \cite[12.13]{BZ}).  Then note that
the reflection $-K$ of $K$ is obtained by replacing $e$ with $r-e$ and
$\beta_i$ with $\alpha_i-\beta_i$.

Also, $K$ is a knot if and only if either exactly one of
$\alpha_1,\dots,\alpha_r$ is even, or if all of
$\alpha_1,\dots,\alpha_r,e+\sum_{i=1}^r\beta_i$ are odd.  We will
restrict our attention to knots.  Since $\delta(K)$ is determined
for alternating knots by Theorem~\ref{alt}, and since
$\delta(-K)=-\delta(K)$, we also restrict to $1\le
e\le\lfloor r/2\rfloor$.

The double cover of $S^3$ branched along $K$ is a Seifert fibred
space which is given as the boundary of a plumbing of disk bundles
over $S^2$ (see for example \cite{OwS}).  This plumbing is
determined (nonuniquely) by the Montesinos invariants which
specify $K$. After possibly reflecting $K$ we may choose the
plumbing so that its intersection pairing is represented by a
negative-definite matrix $Q$.  It then follows from
\cite[Corollary 1.5]{OS3} that
$$\delta(K)=\displaystyle{
\mbox{max}\left\{\frac{w^TQw +\rk(Q)}{2}\right\}},$$ where the
maximum is taken over all $w\in\zz^{\rk(Q)}$ which are
characteristic for $Q$ in the sense of Definition~\ref{char}.

This permits an algorithmic computation of $\delta(K).$ We ran a
{\it Maple} program (partly written by Sa\v{s}o Strle) to find $\delta$ for all Montesinos knots $M(1;
(\alpha_1, \beta_1), (\alpha_2, \beta_2), (\alpha_3, \beta_3))$
with $\alpha_i \leq 7.$ According to the table in \cite{Kaw},
these include all nonalternating Montesinos knots with up to $10$ crossings.

The result is that $\delta = \sigma'$ for all but the knots listed in Table \ref{table:mont}.

\begin{table}
\begin{center}
\begin{tabular}{|c|c|c|c|}
\hline
Montesinos knot & {\it Knotscape} notation & $\delta$ & $\sigma'$ \\
\hline
$M(1; (3,1), (3,1), (4,1))$ & $10n27 = 10_{139}$ & $-1$ & 3 \\
$M(1; (3,1), (3,1), (5,2))$ & $10n14 = 10_{145}$ & 3 & $-1$ \\
$M(1; (2,1), (3,1), (7,1))$ & $12n242$ & 0 & 4 \\
$M(1; (5,2), (5,2), (5,2))$ & $12n276$ & 2 & $-2$ \\
$M(1; (3,1), (4,1), (5,1))$ & $12n472$ & 0 & 4 \\
$M(1; (3,1), (3,1), (6,1))$ & $12n574$ & 0 & 4 \\
$M(1; (2,1), (5,1), (5,1))$ & $12n725$ & 0 & 4 \\
$M(1; (3,2), (5,1), (5,1))$ & $13n3596$ & 4 & 0 \\
$M(1; (5,2), (5,2), (6,1))$ & $14n6349$ & $-2$ & 2 \\
$M(1; (3,1), (4,1), (7,1))$ & $14n12201$ & 1 & 5 \\
$M(1; (3,1), (5,1), (6,1))$ & $14n15856$ & 1 & 5 \\
$M(1; (4,1), (5,1), (5,1))$ & $14n24551$ & 1 & 5 \\
$M(1; (5,1), (5,3), (6,1))$ & $15n74378$ & 0 & 4 \\
$M(1; (5,1), (5,1), (6,1))$ & $16n931575$ & 2 & 6 \\
\hline
\end{tabular}
\vskip5mm
\begin{narrow}{0.3in}{0.3in}
\caption{
{\bf{Montesinos knots for which $\delta\ne\sigma'$.}}
}
\label{table:mont}
\end{narrow}
\end{center}
\end{table}

We used the program {\it Knotscape} \cite{HT} to identify these knots. For the convenience of the
reader, the two ten crossing knots are also shown in Rolfsen's notation \cite{Ro}.

Shumakovitch's {\it KhoHo} package \cite{S} showed that all these knots are not H-thin, in agreement
with Conjecture~\ref{conj}.

\subsection {Torus knots.}
The double branched cover of the torus knot $T(p,q)$ is the Brieskorn sphere $\Sigma(2,p,q).$ This is a Seifert
fibration, so we can use the same formula as in the Montesinos case. Therefore, for every $p$ and $q,$ we have an
algorithmic way of computing $\delta$. For example, the double cover of $T(3,4)$ is $\Sigma(2,3,4),$ 
which gives $\delta(T(3,4)) = \sigma'(T(3,4)) = -3.$

Two infinite classes of knots for which $\delta \neq \sigma'$ are provided by the torus
knots $T(3, 6n+1)$ for $n \geq 0,$ and $T(3, 6n-1)$ for $n \geq 1.$ Indeed, using the recurrence
formulae for the signatures of torus knots in \cite{GLM} we can easily compute $\sigma(T(3, 6n\pm 1))
= 8n.$ Therefore, $\sigma'(T(3, 6n\pm 1)) = -4n.$ On the other hand, the double branched covers are
the Brieskorn spheres $\Sigma (2,3,6n+1),$ whose Heegaard Floer homologies were computed by Ozsv\'ath
and Szab\'o in \cite{OS2}. Their calculations show that $\delta(T(3,6n+1)) = 0$ and $\delta(T(3,
6n-1)) = -4$ for all $n.$

\section {The Ozsv\'ath-Szab\'o technique}

In \cite{OS8}, Ozsv\'ath and Szab\'o managed to compute the correction terms $d$
for $\Sigma(K)$ for four knots of ten crossings that are neither alternating nor Montesinos. Their
method was to find explicit sharp four-manifolds with boundary $\Sigma(K).$ 
If we focus on the spin structure, we see that the results
$$ \delta(10_{148}) = 1, \ \delta(10_{151}) = -1, \ \delta(10_{158}) = 0, \ \delta(10_{162}) = 1 $$
coincide with $\sigma'$ in all four examples.

A rational homology sphere $Y$ is called an $L$-space if its reduced Heegaard Floer homology
$HF_{\mathrm{red}}(Y)$ vanishes, so that its Heegaard Floer invariants resemble those of a lens space.
A negative-definite four-manifold $X$ is called sharp if its boundary $\partial X=Y$ is an $L$-space, and if
$$d(Y,\tt)=\displaystyle
\mbox{max}\left\{\left.\frac{c_1(\ss)^2+b_2(X)}{4}\ \right|\ \ss\in\spc(X),\ss|_{Y}=\tt\right\}$$
for all $\tt\in\spc(Y)$.  This means that the correction terms of $Y$ can be computed from
the intersection pairing of $X$.  We say a negative-definite form $Q$ is sharp for $Y$ if $Q$ is
the intersection pairing of a sharp four-manifold bounded by $Y$.

The following proposition summarizes the above-mentioned technique. The proof of this proposition
is due to Ozsv\'ath and Szab\'o and may be found in \cite[\S 7.2]{OS8}.

\begin{proposition}
Let $k\ge2$.
Suppose that a knot $K$ has a nonalternating projection which contains a subdiagram which is a $-k/2$ twist
(a $2$-strand braid with $k$ crossings, as shown below) 
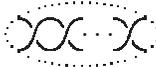
\begin{center}
\ifpic \leavevmode \xygraph{ !{0;/r1.0pc/:}
%%% first diagram: -k/2 twist
!{\htwistneg}!{\htwistneg}[rd(0.5)]!{*{\dots}}[ru(0.5)]!{\htwistneg}[l(2)d(0.5)]!{\ellipse(2.4,1){.}}
}
\else \vskip 5cm \fi
\end{center}
and that the complement of this subdiagram is alternating.
Let $K',K''$ respectively be the links obtained by replacing the $-k/2$ twist by
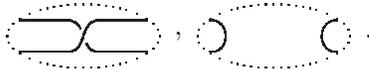
\begin{center}
\ifpic \leavevmode \xygraph{ !{0;/r1.0pc/:}
%%% second diagram: -1/2 twist
[u(0.5)r(4)]
!{\xcaph[1.5]@(0)}[ld]
!{\xcaph[1.5]@(0)}[u(1)r(0.5)]
!{\htwistneg}
!{\xcaph[1.5]@(0)}[ld]
!{\xcaph[1.5]@(0)}[l(1.5)u(0.5)]
!{\ellipse(2.4,1){.}}
%%% third diagram: 1 resolution
[r(3)]!{*{,}}[u(0.5)r]
!{\hcap}[r(4)]!{\hcap[-]}[l(2)d(0.5)]!{\ellipse(2.4,1){.}}
[r(3)]!{*{.}}
}
\else \vskip 5cm \fi
\end{center}
Note that the resulting projection of $K''$ is alternating.  Let $G$ be the negative-definite
$m\times m$ Goeritz matrix associated to $K''$ as in Section \ref{sec:alt}, where we take
$v_0$ and $v_m$ to be the vertices in the white graph associated to the regions at each end
of the subdiagram shown above.  Let $G_r$ be the matrix
$$\left(\begin{array}{cccl}
&&\hfill\vline&\\
&G&\hfill\vline&\\
&&\hfill\vline&1\\
\hline
&&1\,\,\hfill\vline& r\\
\end{array}\right).$$
Then $G_{-1}$ is the intersection pairing of a four-manifold bounded by $\Sigma(K')$, and
$G_{-k}$ is the intersection pairing of a manifold bounded by $\Sigma(K)$.
If $G_{-1}$ is negative-definite and sharp for $\Sigma(K')$ then $G_{-k}$ is sharp for $\Sigma(K)$.
\end{proposition}

\begin {figure}
\begin {center}
\begin{picture}(0,0)%
\includegraphics{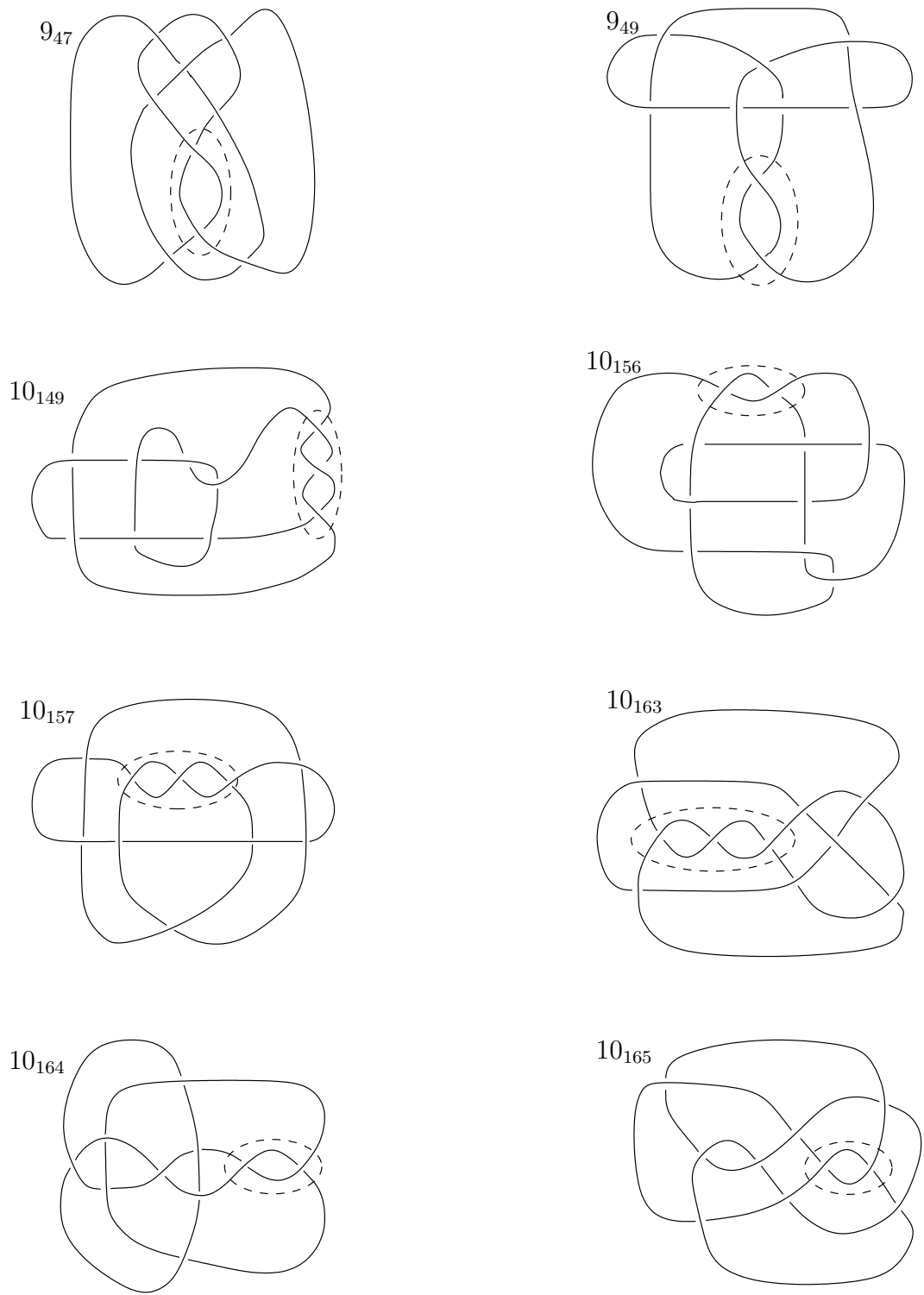}%
\end{picture}%
\setlength{\unitlength}{3947sp}%
\begingroup\makeatletter\ifx\SetFigFont\undefined%
\gdef\SetFigFont#1#2#3#4#5{%
  \reset@font\fontsize{#1}{#2pt}%
  \fontfamily{#3}\fontseries{#4}\fontshape{#5}%
  \selectfont}%
\fi\endgroup%
\begin{picture}(6655,9378)(1801,-9352)
\put(1801,-2836){\makebox(0,0)[lb]{\smash{{\SetFigFont{12}{14.4}{\rmdefault}{\mddefault}{\updefault}{\color[rgb]{0,0,0}$10_{149}$}%
}}}}
\put(1876,-5161){\makebox(0,0)[lb]{\smash{{\SetFigFont{12}{14.4}{\rmdefault}{\mddefault}{\updefault}{\color[rgb]{0,0,0}$10_{157}$}%
}}}}
\put(1801,-7711){\makebox(0,0)[lb]{\smash{{\SetFigFont{12}{14.4}{\rmdefault}{\mddefault}{\updefault}{\color[rgb]{0,0,0}$10_{164}$}%
}}}}
\put(6076,-7636){\makebox(0,0)[lb]{\smash{{\SetFigFont{12}{14.4}{\rmdefault}{\mddefault}{\updefault}{\color[rgb]{0,0,0}$10_{165}$}%
}}}}
\put(6001,-2611){\makebox(0,0)[lb]{\smash{{\SetFigFont{12}{14.4}{\rmdefault}{\mddefault}{\updefault}{\color[rgb]{0,0,0}$10_{156}$}%
}}}}
\put(2026,-211){\makebox(0,0)[lb]{\smash{{\SetFigFont{12}{14.4}{\rmdefault}{\mddefault}{\updefault}{\color[rgb]{0,0,0}$9_{47}$}%
}}}}
\put(6151,-136){\makebox(0,0)[lb]{\smash{{\SetFigFont{12}{14.4}{\rmdefault}{\mddefault}{\updefault}{\color[rgb]{0,0,0}$9_{49}$}%
}}}}
\put(6151,-5086){\makebox(0,0)[lb]{\smash{{\SetFigFont{12}{14.4}{\rmdefault}{\mddefault}{\updefault}{\color[rgb]{0,0,0}$10_{163}$}%
}}}}
\end{picture}%

\caption {The relevant $-1$ or $-3/2$ twist is indicated on each knot.}
\label {fig:knots}
\end {center}
\end {figure}

For each knot $K$ in Figure \ref{fig:knots} we find that $\Sigma(K')$ is a lens space whose correction
terms may be computed using the formula from \cite{OS2}, so that the condition on $G_{-1}$ may be
checked.  Then the matrix $G_{-k}$ may be used to compute the correction terms of $\Sigma(K)$.
In particular we find that $\delta=\sigma'$ for all of these knots.

Putting this together with the calculations for alternating, Montesinos, and torus knots, this covers all 
but eight of the knots up to ten crossings, and establishes the result in Theorem~\ref{nine}. 

\medskip 

In Table \ref{table:smallknots} we exhibit all the knots of up to ten crossings for which either some of the four 
invariants $\delta, \sigma', \tau$ and $s'$ are different, or some of their values are unknown. In other 
words, for all the knots up to ten crossings that are not shown, we know that $\delta = \sigma' = \tau = 
s'.$

The notation for knots is the one in \cite{Lt}, the same as Rolfsen's but with the last five knots 
translated to account for the Perko pair. The values for $\sigma', \tau$ and $s'$ are also taken from \cite{Lt}.

\begin{table}[!ht]
\begin{center}
\begin{tabular}{|c|c|c|c|c|c|}
\hline
Knot & Type & $\delta$ & $\sigma'$ & $\tau$ & $s'$\\
\hline
$ 9_{42} $ 	& H-thick & $-1$ & $-1$ & 0 & 0 \\
$ 10_{132}$	& H-thick & 0 & 0 & $-1$ & $-1$ \\
$ 10_{136}$	& H-thick & $-1$ & $-1$ & 0 & 0\\
$ 10_{139}$	& H-thick & $-1$ & 3 & 4 & 4 \\
$ 10_{141}$	& H-thick & 0 & 0 & ? & 0 \\
$ 10_{145}$	& H-thick & 3 & $-1$ & $-2$ & $-2$ \\
$ 10_{150}$	& H-thin & ? & 2 & 2 & 2 \\
$ 10_{152}$	& H-thick & ? & 3 & 4 & 4 \\
$ 10_{153}$	& H-thick & ? & 0 & 0 & 0 \\
$ 10_{154}$	& H-thick & ? & 2 & 3 & 3 \\
$ 10_{155}$	& H-thin & ? & 0 & 0 & 0 \\
$ 10_{159}$	& H-thin & ? & 1 & 1 & 1 \\
$ 10_{160}$	& H-thin & ? & 2 & 2 & 2 \\
$ 10_{161}$	& H-thick & ? & 2 & 3 & 3 \\
\hline
\end{tabular}
\vskip5mm
\begin{narrow}{0.3in}{0.3in}
\caption{
{\bf{Knots with 10 crossings or less for which either $\delta$, $\sigma'$, $\tau$, $s'$ are not all equal, 
or some of their values are not known.}}
}
\label{table:smallknots}
\end{narrow}
\end{center}
\end{table}

\section {Whitehead doubles}

The double branched coverings of doubles of knots were studied extensively in \cite{MW}. In
particular, it is interesting to note that the double branched cover of a twisted double of a prime
knot $K$ determines $K$ among all prime knots.

\begin{figure}[htbp]
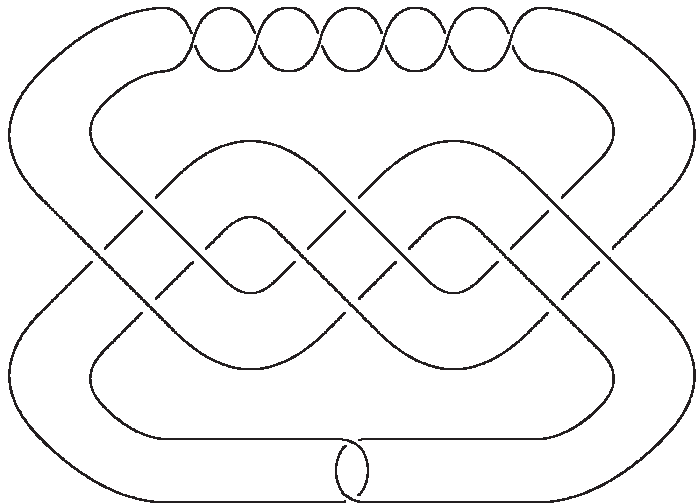
  %%% Wh(trefoil)
\begin{center}
\ifpic \leavevmode \xygraph{ !{0;/r2pc/:}
%%% 6 crossings in middle, top
!{\htwistneg} !{\htwistneg} !{\htwistneg} !{\htwistneg}
!{\htwistneg} !{\htwistneg}
%%% top left
[l(7.9)]!{\xbendl[-2]} [r]!{\xbendl-} [ul]!{\xcapv[-2]}
[r(1)u(0.5)]!{\xcapv[-0.9]}
%%% 1st double positive crossing
[d(0.2)r(0.3)]!{\xunderv[0.8]} [l(0.8)u(0.2)]!{\xunderv[0.8]}
[u(1)r(1.6)]!{\xunderv[0.8]} [l(0.8)u(0.2)]!{\xunderv[0.8]}
%%% curved pieces: 2 large
[u(2.6)r(0.8)]!{\xcaph[2.4]} [l(1)d(2.4)]!{\xcaph[-2.4]}
%%% 2nd double positive crossing
[u(2.4)r(1.4)]!{\xunderv[0.8]} [l(0.8)u(0.2)]!{\xunderv[0.8]}
[u(1)r(1.6)]!{\xunderv[0.8]} [l(0.8)u(0.2)]!{\xunderv[0.8]}
%%% curved pieces: 2 large
[u(2.6)r(0.8)]!{\xcaph[2.4]} [l(1)d(2.4)]!{\xcaph[-2.4]}
%%% 3rd double positive crossing
[u(2.4)r(1.4)]!{\xunderv[0.8]} [l(0.8)u(0.2)]!{\xunderv[0.8]}
[u(1)r(1.6)]!{\xunderv[0.8]} [l(0.8)u(0.2)]!{\xunderv[0.8]}
%%% top right
[u(5.3)r(0.1)]!{\xbendr[-2]} [l]!{\xbendr-} [ur]!{\xcapv[2]}
[l(1)u(0.5)]!{\xcapv[0.9]}
%%% bottom right
[d(3.3)l(1)]!{\xbendl[2]} [u]!{\xbendl} [ur(2)]!{\xcapv[2]}
[l(1)u(0.4)]!{\xcapv[0.9]}
%%% bottom left
[u(0.6)l(8.8)]!{\xbendr[2]} [u]!{\xbendr} [ul(2)]!{\xcapv[-2]}
[r(1)u(0.4)]!{\xcapv[-0.9]}
%%% bottom, with a positive clasp
[r(1)d(0.4)]!{\xcaph[2.7]@(0)} [dl]!{\xcaph[2.7]@(0)}
[u(1)r(1.65)]!{\vover[0.5]} [u(0.5)]!{\vover[-0.5]}
[u(1.5)r(0.5)]!{\xcaph[2.65]@(0)} [dl]!{\xcaph[2.65]@(0)}
%%% now fill in some gaps: first small curves
[u(4.3)l(0.05)]!{\xcaph[0.8]} [l(4.2)]!{\xcaph[0.8]}
[l(1)d(0.8)]!{\xcaph[-0.8]} [r(2.2)]!{\xcaph[-0.8]}
%%% finally some missing segments
[u(0.8)r(2.2)]!{\sbendh[0.5]@(0)}
[d(0.2)l(1.8)]!{\sbendh[0.3]@(0)}
[d(2.6)l(0.2)]!{\sbendv[0.5]@(0)}
[u(0.2)l(1.8)]!{\sbendv[0.3]@(0)}
[u(3.7)l(8.5)]!{\sbendv[0.3]@(0)} [u(0.4)l(2)]!{\sbendv[0.5]@(0)}
[d(1.4)l(0)]!{\sbendh[0.3]@(0)} [d(0.4)l(2)]!{\sbendh[0.5]@(0)} }
\else \vskip 5cm \fi
\begin{narrow}{0.3in}{0.3in}
\caption{ {\bf The untwisted Whitehead double of the right handed
trefoil, with a positive clasp.} } \label{fig:whtref}
\end{narrow}
\end{center}
\end{figure}

In this section we focus on the untwisted Whitehead double $Wh(K)
= Wh^+(K)$ of a knot $K$ with a positive clasp.  The untwisted
double of the right handed trefoil is shown in Figure
\ref{fig:whtref}.  (The three negative twists shown in the diagram are
a consequence of the fact that the chosen diagram of the trefoil
has writhe three.)
Note that the untwisted double $Wh^-(K)$ with a
negative clasp is simply $-Wh^+(-K),$ so our discussion can also
be phrased in terms of $Wh^-.$

Denote by $K^r$ the knot $K$ with its string orientation reversed. 
The following surgery description of $\Sigma(Wh(K))$ is 
useful:

\begin {proposition}
\label {surg} $\Sigma(Wh(K))$ is the manifold obtained by Dehn
surgery on the knot $K \# K^r$ in $S^3$ with framing $1/2.$
\end {proposition}

\begin{figure}[htbp]
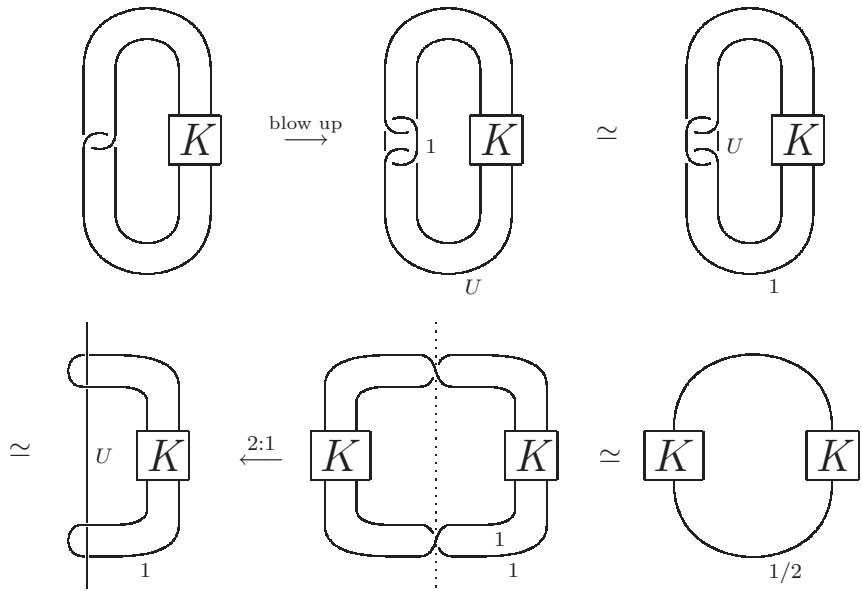
 %%% double cover of Wh(K)
\begin{center}
\ifpic \leavevmode \xygraph{ !{0;/r1.0pc/:}
%%%%%%%%%%%%%%%% first draw Wh(K)
%%top of Wh(K)
!{\vcap[4]}[r] !{\vcap[2]} % two maxima
%%left half of Wh(K)
[l] !{\xcapv[2]@(0)} [ur] !{\xcapv[2]@(0)} % 2 vertical segments
[ld]!{\hover[0.5]}[l(0.5)]!{\hover[-0.5]} % positive clasp
[l(1.5)d(0.5)]!{\xcapv[2]@(0)} [ur] !{\xcapv[2]@(0)} % 2 vertical segments
%%bottom of Wh(K)
[dl]!{\vcap[-4]} [r] !{\vcap[-2]} % two minima
%%right half of Wh(K)
[r(2)u(4.5)] !{\xcapv[1.45]@(0)} [ur] !{\xcapv[1.45]@(0)} % 2 vertical segments
[l(0.5)d(1.2)] !{*+[F]{\mbox{\it\LARGE K}}} % K in box
[l(0.5)d(0.8)] !{\xcapv[1.5]@(0)} [ur] !{\xcapv[1.5]@(0)} % 2 vertical segments
%%blow up arrow
[r(3)u(2.1)] !{*{\stackrel{\mbox{\tiny blow
up}}{\longrightarrow}}}
%%%%%%%%%%%%%%%% unknotted Wh(K)
%%top
[u(1.9)r(2.5)]!{\vcap[4]}[r] !{\vcap[2]} % two maxima
%%left half
[l] !{\xcapv[1.5]@(0)} [ur] !{\xcapv[1.5]@(0)} % 2 vertical segments
[l(1)d(0.5)]!{\hunder[0.5]}[l(0.5)]!{\hunder[-0.5]} % clasp
[l(1.5)d(0.5)]!{\xcapv[.5]@(0)} [ur] !{\xcapv[.5]@(0)<{1}} % 2 vertical segments
[l(1)u(0.5)]!{\hover[0.5]}[l(0.5)]!{\hover[-0.5]} % clasp
[l(1.5)d(0.5)]!{\xcapv[1.5]@(0)} [ur] !{\xcapv[1.5]@(0)} % 2 vertical segments
%%bottom
[d(0.5)l]!{\vcap[-4]>{U}} [r] !{\vcap[-2]} % two minima
%%right half
[r(2)u(4.5)] !{\xcapv[1.45]@(0)} [ur] !{\xcapv[1.45]@(0)} % 2 vertical segments
[l(0.5)d(1.2)] !{*+[F]{\mbox{\it\LARGE K}}} % K in box
[l(0.5)d(0.8)] !{\xcapv[1.5]@(0)} [ur] !{\xcapv[1.5]@(0)} % 2 vertical segments
%%isotopic
[r(3)u(2)] !{*{\simeq}}
%%%%%%%%%%%%%% unknotted Wh(K) again (switch unknots)
%%top
[u(2)r(2.5)]!{\vcap[4]}[r] !{\vcap[2]} % two maxima
%%left half
[l] !{\xcapv[1.5]@(0)} [ur] !{\xcapv[1.5]@(0)} % 2 vertical segments
[l(1)d(0.5)]!{\hover[0.5]}[l(0.5)]!{\hover[-0.5]} % clasp
[l(1.5)d(0.5)]!{\xcapv[.5]@(0)} [ur] !{\xcapv[.5]@(0)<{U}} % 2 vertical segments
[l(1)u(0.5)]!{\hunder[0.5]}[l(0.5)]!{\hunder[-0.5]} % clasp
[l(1.5)d(0.5)]!{\xcapv[1.5]@(0)} [ur] !{\xcapv[1.5]@(0)} % 2 vertical segments
%%bottom
[d(0.5)l]!{\vcap[-4]>{1}} [r] !{\vcap[-2]} % two minima
%%right half
[r(2)u(4.5)] !{\xcapv[1.45]@(0)} [ur] !{\xcapv[1.45]@(0)} % 2 vertical segments
[l(0.5)d(1.2)] !{*+[F]{\mbox{\it\LARGE K}}} % K in box
[l(0.5)d(0.8)] !{\xcapv[1.5]@(0)} [ur] !{\xcapv[1.5]@(0)} % 2 vertical segments
%%isotopic
[l(25)d(8)] !{*{\simeq}}
%%%%%%%%%%%%%% unknotted Wh(K) again, U vertical
%%clasp around vertical U
[u(3)r(2)]!{\hcap[-1]} % opening to right
[r(0.2)]!{\xcaph[1]@(0)} % horizontal segment
[d(1)l(1.2)]!{\xcaph[1.2]@(0)} % horizontal segment
%%curve at top right
[u]!{\xbendr[-2]@(0)}[l] !{\xbendr-@(0)} [ur]!{\xcapv[0.5]@(0)}
%%now K in the middle
[u(0.5)l]!{\xcapv[0.9]@(0)} [ur] !{\xcapv[0.9]@(0)} % 2 vertical segments
[l(0.5)d(0.65)] !{*+[F]{\mbox{\it\LARGE K}}} % K in box
[l(0.5)d(0.8)] !{\xcapv[0.9]@(0)} [ur] !{\xcapv[0.9]@(0)} % 2 vertical segments
%%curve at bottom right
[u(0.1)]!{\xcapv[0.5]@(0)} [l(2)u(1.5)] !{\xbendl@(0)}
[u]!{\xbendl[2]@(0)>{1}}
%%clasp around vertical U
[l(1)]!{\hcap[-1]} % opening to right
[r(0.2)d(1)]!{\xcaph[1]@(0)} % horizontal segment
[u(1)l(1.2)]!{\xcaph[1.2]@(0)} % horizontal segment
%%now draw U
[u(6.35)l(0.9)]!{\xcapv[1.9]@(0)} % vertical segment
[d(1.1)]!{\xcapv[4.15]@(0)|{U}} % vertical segment
[d(3.35)]!{\xcapv[1.9]@(0)}
%%%%%%%%%%%%%%%%% now draw double cover
[u(3.45)r(5.5)]!{*{\stackrel{2:1}{\longleftarrow}}}
%%clasp becomes crossing
[u(3)r(5)]!{\htwist} % crossing top middle
[d(1)]!{\xcaph[1]@(0)} % horizontal segment
[u(1)l(1)]!{\xcaph[1]@(0)} % horizontal segment
%%curve at top right
!{\xbendr[-2]@(0)}[l] !{\xbendr-@(0)} [ur]!{\xcapv[0.5]@(0)}
%%now K in the middle
[u(0.5)l]!{\xcapv[0.9]@(0)} [ur] !{\xcapv[0.9]@(0)} % 2 vertical segments
[l(0.5)d(0.65)] !{*+[F]{\mbox{\it\LARGE K}}} % K in box
[l(0.5)d(0.8)] !{\xcapv[0.9]@(0)} [ur] !{\xcapv[0.9]@(0)} % 2 vertical segments
%%curve at bottom right
[u(0.1)]!{\xcapv[0.5]@(0)} [l(2)u(1.5)] !{\xbendl@(0)>{1}}
[u]!{\xbendl[2]@(0)>{1}}
%%clasp becomes crossing
[l(2)]!{\htwistneg} % crossing bottom middle
[d(1)]!{\xcaph[1]@(0)} % horizontal segment
[u(1)l(1)]!{\xcaph[1]@(0)} % horizontal segment
%%now draw U
[u(6.35)l(1.5)]!{\knotstyle{.}}!{\xcapv[8.3]@(0)}!{\knotstyle{-}} % vertical dotted segment
%%left side
[d(1)l(1.5)]!{\xcaph[1]@(0)} % horizontal segment
[u(1)l(1)]!{\xcaph[1]@(0)} % horizontal segment
%%curve at top left
[l(3)]!{\xbendl[-2]@(0)} !{\xcapv[0.5]@(0)} [r(1)u(1)]
!{\xbendl-@(0)}
%%now K in the middle
[u(0.5)l]!{\xcapv[0.9]@(0)} [ur] !{\xcapv[0.9]@(0)} % 2 vertical segments
[l(0.5)d(0.65)] !{*+[F]{\mbox{\it\LARGE $K$}}} % K in box
[l(0.5)d(0.8)] !{\xcapv[0.9]@(0)} [ur] !{\xcapv[0.9]@(0)} % 2 vertical segments
%%curve at bottom right
[u(0.6)l(1)]!{\xbendr[2]@(0)} [u]!{\xbendr@(0)}
[l(2)u(0.5)]!{\xcapv[0.5]@(0)}
[d(0.5)r(2)]!{\xcaph[1]@(0)} % horizontal segment
[u(1)l(1)]!{\xcaph[1]@(0)} % horizontal segment
%%isotopic
[r(6)u(2.15)] !{*{\simeq}}
%%%%%%%%%%%%%% final Dehn surgery diagram
[u(0.8)r(2)]!{\vcap[5]} % large maximum
[d(0.75)]!{*+[F]{\mbox{\it\LARGE $K$}}} [r(5)] !{*+[F]{\mbox{\it\LARGE K}}} % K in boxes
[l(5)d(0.8)]!{\vcap[-5]>{1/2}} % large minimum
} \else \vskip 5cm \fi
\begin{narrow}{0.4in}{0.2in}
\caption{ {\bf{Proof that the double branched cover of $Wh(K)$ is
$(K\#K^r)_{1/2}$.}}
Note that in the last diagram $K$ is traversed twice with different orientations, so that
the knot shown is $K\#K^r$ rather than $K\#K$.} 

\label{2cover}
\end{narrow}
\end{center}
\end{figure}

\noindent {\bf Proof.}  The knot $Wh(K)$ may be unknotted by
changing one crossing, so as to undo the clasp.  A theorem
originally due to Montesinos tells us that $\Sigma(Wh(K))$ is
given by Dehn surgery on a knot with half-integral framing. An
algorithm is described in \cite[Lemma 3.1]{Ow} to obtain a Dehn
surgery diagram for the double branched cover of any knot $C$ by
first obtaining a Kirby calculus diagram of a disk bounded by $C$
in a blown-up four-ball.  The description of $\Sigma(Wh(K))$
follows from the application of this algorithm to $Wh(K)$. (This
is illustrated in Figure~\ref{2cover}.) $\hfill \fin$

\medskip

In \cite{R2}, Rasmussen introduced a set of invariants $h_i$ for a knot $K \subset S^3,$ where $i \in
\zz.$ They are nonnegative integers that encode the information in the maps in the surgery exact
triangles for $K.$ According to \cite[Section 7.2]{R2} and \cite[Section 2.2]{R4}, the correction
terms $d$ for manifolds obtained by integral surgery on a knot $K \subset S^3$ can be computed in
terms of the $h_i$'s. In particular, the $d$ invariant of the $-1$ surgery on $K$ is equal to
$2h_0.$

\begin {proposition}
\label {dd}
Let $K_{-1/2}$ be the manifold obtained from $S^3$ by surgery on the knot $K$ with framing $-1/2.$
Then $d(K_{-1/2}) = 2h_0(K).$ (Note that $K_{-1/2}$ has a unique $\spc$ structure, which we drop
from the notation.)
\end {proposition}

\noindent {\bf Proof.} Let $K_0$ and $K_{-1}$ be the results of $0$ and $-1$ surgery on $K \subset S^3,$
respectively. Consider the exact triangle with twisted coefficients \cite{OS1}:
\begin {equation}
\label {equa1}
\begin {CD}
{\underline {HF}}^+(K_0) \to HF^+(S^3)[T, T^{-1}]
@>{F_{W}}>> HF^+(K_{-1})[T, T^{-1}] @>{F_{W_0}}>> {\underline {HF}}^+(K_0)
\end {CD}
\end {equation}

We tensor everything with $\qq$ for simplicity, so that (\ref{equa1}) is a long exact sequence
of modules over the ring $\qq[U, T, T^{-1}].$ The map $F_{W}$ is induced by the cobordism $W$ corresponding
to the $-1$ surgery on $K.$  Choose a generator of $H^2(W)$ and denote by $\tt_i$ the $\spc$
structure on $W$ with $c_1(\tt_i) = 2i -1.$ There are maps $F_{W, \tt_i}: HF^+(S^3) \to
HF^+(K_{-1})$ associated to each $\spc$ structure, and
$$ F_W(x) = \sum T^i \cdot F_{W, \tt_i} (x).$$

Since $S^3$ and $K_{-1}$ are integral homology spheres, $F_{W, \tt_i}$ are graded maps. More
precisely, they shift the absolute grading on $HF^+$ by $(c_1(\tt_i)^2 + 1)/4 = i - i^2.$ In
particular, the maximal degree shift is zero.

For every integral homology sphere $Y$, the Heegaard Floer homology $HF^+(Y; \qq)$ is a
$\qq[U]$-module that decomposes as $\TT^+_{d(Y)} \oplus HF_{\mathrm{red}}(Y).$ Here $HF_{\mathrm{red}}$ is finite
dimensional as a $\qq$-vector space, and $\TT^+_{d(Y)} \cong \qq[U, U^{-1}]/\qq[U]$ with $U$ reducing
the degree by $2$ and the bottom-most element in $\TT^+_{d(Y)}$ being in absolute degree $d(Y).$ The
piece $\TT^+_{d(Y)}$ can be described as the image of $HF^{\infty}(Y)$ in $HF^+(Y),$ or as the image of
$U^k$ in $HF^+(Y)$ for $k$ sufficiently large.

By $U$-equivariance, the map $F_W$ must take $HF^+(S^3)[T, T^{-1}] = \TT^+_0[T, T^{-1}]$ into
$\TT^+_{d(K_{-1})}$ $[T, T^{-1}].$ The
number $h_0 = h_0(K)$ is the rank of the kernel of $F_W$ as a $\qq[T, T^{-1}]$ module, which is the span of
$\{1, U^{-1}, \dots , U^{-h_0+1} \}.$ Let $x$ be a generator of the part of $\TT^+_{d(K_{-1})}$ sitting in
the bottom-most degree $d(K_{-1}) = 2h_0.$ For all $k$ we have $$ F_W(U^{-h_0-k}) = P(T)U^{-k}x + \text{
lower degree terms}, $$
where $P(T)$ is a fixed polynomial in $T.$ Since the cobordism $W$ is
negative definite, each of the $\spc$ structures $\tt_i$ (in particular $\tt_0$ and $\tt_1$) induces an
isomorphism on $HF^{\infty}.$ It follows that $P(T) = aT+b$ with $a, b \neq 0.$

The conclusion is that for any element $z$ in the image of $F_W$ (which is also the
kernel of $F_{W_0}$), the part of $z$ sitting in highest degree is always a multiple of $P(T).$
Note that $P(T)$ is not invertible in $\qq[T, T^{-1}].$

Now let $K' \subset K_{-1}$ be the core of the solid torus glued in for the $-1$ surgery. Then $-1$
surgery on $K' \subset K_{-1}$ produces the manifold $K_{-1/2},$ while zero surgery on $K' \subset
K_{-1}$ gives $K_0$ again.

We consider the exact triangle with twisted coefficients for the $-1$ surgery on $K' \subset
K_{-1}:$
\begin {equation}
\label {equa2}
\begin {CD}
\dots \to {\underline {HF}}^+(K_0) @>{F_{Z}}>> HF^+(K_{-1})[T, T^{-1}]
@>{F_{W'}}>> HF^+(K_{-1/2})[T, T^{-1}] \to \dots
\end {CD}
\end {equation}

The map $F_{Z}$ is induced by a cobordism $Z$ from $K_0$ to $K_{-1}.$ Drawing inspiration from the
arguments in \cite[Lemma 2.9]{OSe} and \cite[Lemma 4.5]{OS8}, we compose $Z$ with the cobordism $W_0$ from
$K_{-1}$ to $K_0$ that appeared in (\ref{equa1}). Then $W_0 \circ Z$ has an alternative factorization as
the two-handle addition from $K_0$ to $K_0 \# (S^2 \times S^1),$ followed by another cobordism in which
the generator of $H_1(S^2 \times S^1)$ becomes null-homologous. This implies that $F_{W_0} \circ F_Z = F_{W_0
\cup Z} = 0.$

It follows from here that the image of $F_Z$ (which is also the kernel of $F_{W'}$) lies in the kernel of
$F_{W_0}.$ Hence every element in the kernel of $F_{W'}$ must have its highest degree part a multiple of
$P(T).$ In particular, since the bottom-most generator $x$ of $\TT^+_{d(K_{-1})}$ is not of this form, we
must have $F_{W'}(x) \neq 0.$

Note that $K_{-1}$ and $K_{-1/2}$ are integral homology spheres, and $F_{W'}$ is a sum of graded maps as
before, with the maximal shift in absolute grading being zero. The map $F_{W'}$ must take
$\TT^+_{d(K_{-1})}[T, T^{-1}]$ to $\TT^+_{d(K_{-1/2})}[T, T^{-1}].$  Also, $W'$ is a negative definite
cobordism, which means that $F_{W'}(x)$ must have its highest degree part in grading $2h_0 = $ degree$(x).$
On the other hand, we have $U\cdot F_{W'}(x) = F_{W'}(Ux) = 0,$ hence $F_{W'}(x)$ lives in the bottom-most
degree $d(K_{-1/2}).$ We conclude that $d(K_{-1/2}) = 2h_0. \hfill \fin $

\medskip

Applying this to the mirror image of $K$ we get:

\begin {corollary}
\label {de}
Let $K_{1/2}$ be the manifold obtained from $S^3$ by surgery on the knot $K$ with framing $1/2.$
Then $d(K_{1/2}) = -2h_0(-K).$ 
\end {corollary}

\medskip

\noindent {\bf Proof of Theorem~\ref{double}.} Proposition~\ref{surg} and Corollary~\ref{de} imply that
$\delta(Wh(K)) = -4h_0(-(K \# K^r)).$ Since $h_0$ is always nonnegative, we must have $\delta(Wh(K)) \leq 0.$
Proposition 7.7 in \cite{R2} says that $h_0(K) > 0$ whenever $\tau(K) < 0.$ Applying this to $-(K \# K^r)$
and using the fact that $\tau(-(K \# K^r)) = -2\tau(K),$ we get that $\delta(Wh(K)) < 0$ for $\tau(K) > 0.$

According to \cite{OS5}, alternating knots are perfect in the sense of \cite[Definition 6.1]{R2}, i.e.
their knot Floer homology is supported on one diagonal. For a perfect knot $K, \ h_0(K) = \max
\{\lceil -\tau(K)/2 \rceil, 0 \}$ by \cite[Corollary 7.2]{R2}. The property of being perfect is
preserved under taking connected sums by \cite[Corollary 6.2]{R2}. Also, knot Floer homology is insensitive to 
orientation reversal of the knot. Therefore, if $K$ is alternating
(or, more generally, perfect) then $-(K \# K^r)$ is perfect and $\delta(Wh(K)) = -4\max \{\tau(K), 0 \}.
\hfill \fin$

\medskip

\noindent {\bf Proof of Corollary~\ref{coronew}.} The knot $K_1=T(2, 2m+1)$ is alternating and has $\tau(K_1) 
= m > 0.$ By Theorem~\ref{double} we have $\delta(Wh(K_1)) = -4m.$ On the other hand, by writing down 
a Seifert matrix it is easy to see that the signature of the Whitehead double of any knot is zero, so 
$\sigma'(Wh(K_1)) =0.$ The results of Livingston and Naik from \cite{LN}, together with the fact that the 
Thurston-Bennequin number $TB(K_1) = 2m-1$ is positive (see for example \cite{N}), imply that $\tau(Wh(K_1)) 
= s'(Wh(K_1)) = 1.$ The same results hold for $K_2 = T(2, 2n+1),$ and the corollary follows from the 
additivity properties of the four invariants. $\hfill \fin$

\medskip 

The fact that $\delta(Wh(T(2, 2m+1)) = -4m$ shows that $\delta(Wh(K))$ can be any nonpositive multiple of 
four. In contrast to this, $\sigma'(Wh(K)) = 0$ for all $K$ and $s'(Wh(K)), \tau(Wh(K)) \in \{ 0, \pm 1\}$ 
because their absolute values are bounded above by the unknotting number of $Wh(K).$

\medskip

\noindent {\bf Proof of Corollary~\ref{coro2}.} Note that every topologically slice knot has $\sigma = 0,$ and 
therefore $\delta \equiv 0 \pmod 4$ by (\ref{congruence}). We claim that the homomorphism $\phi = (\tau, \delta/4): \ccc_{ts} 
\to \zz \oplus \zz$ is surjective. Indeed, we have $\phi(Wh(T(2,3)) = (1, -1)$ and $\phi(Wh(T(2,5)) = (1, -2),$ and these 
values span $\zz \oplus \zz.\ \hfill \fin$

\end{document}